\documentclass{amsproc}

\usepackage{amssymb}
\usepackage{graphicx}

\usepackage{}

\newtheorem{theorem}{Theorem}[section]

\newtheorem{prob}[theorem]{Problem}

\newtheorem{conj}[theorem]{Conjecture}

\theoremstyle{definition}

\theoremstyle{remark}

 \newcommand{\res}{\mbox{\rm Res}}

\newcommand{\C}{\mathbb{C}}
\newcommand{\Z}{\mathbb{Z}}

\newcommand{\one}{\mathbf{1}}

\numberwithin{equation}{section}

\begin{document}

\title[Some open problems]{Some open problems in mathematical 
two-dimensional conformal field theory}


\author{Yi-Zhi Huang}
\address{Department of Mathematics,
Rutgers University,
110 Frelinghuysen Road,
Piscataway, NJ 08854,
USA}
\curraddr{}
\email{yzhuang@math.rutgers.edu}
\dedicatory{Dedicated to Jim Lepowsky and Robert Wilson for their 70th birthday}
\thanks{}


\date{}

\begin{abstract}
We discuss some open problems in a program of constructing and studying 
two-dimensional conformal field theories using the representation theory of vertex operator
algebras. 
\end{abstract}

\maketitle

\section{Introduction}

Quantum field theory has become an active research area in mathematics in the last forty years.
Among all the quantum field theories, topological quantum field theory is the most successful 
in mathematics mainly because the state space of a topological quantum field theory is typically
finite dimensional. Compared with topological quantum field theory,  nontopological quantum field theories and the 
deep mathematical conjectures derived from these theories are much more difficult and are
still quite distant from a complete 
mathematical understanding. 

One of the most famous but also one of the most difficult problems on nontopological quantum field theories is the existence
of four-dimensional quantum Yang-Mills theory and the mass gap problem. On the other hand, two-dimensional conformal 
field theory as the best understood 
nontopological quantum field theory has in fact been greatly developed and 
has also  directly provided ideas and methods for the successful solutions of 
mathematical conjectures and problems. The study of two-dimensional conformal field theory will certainly also
shed light on the other more difficult nontopological quantum field theories such as the four-dimensional Yang-Mills theory.

In a program of constructing and studying 
two-dimensional conformal field theories using the representation theory of vertex operator
algebras, the mathematical foundation of two-dimensional 
conformal field theory has been essentially established. 
In this note, we discuss some open problems in this program.  Some of the problems are well known. 
The others are what the author is interested in and 
believes to be important. 
Most of the problems here are enhanced and modified versions of some problems 
in the slides that the author prepared for the problem session of this conference.

These problems emphasize mostly the general theory rather than special examples.
For many years, compared with the study of special examples, the development of the general theory in this area
has not received the attention that it should have received.  Though it is important to verify that special 
examples satisfy suitable finiteness and/or  complete reducibility conditions,
it is at least equally important to give a construction of conformal field theories,
prove various conjectures proposed by physicists and mathematicians
and to solve open problems that will lead to future development of the theory and 
future interaction with other branches of mathematics and physics. In fact, many results in the general theory including
the solutions of longstanding open problems, for example,
the operator product expansion of (logarithmic) intertwining operators \cite{Hu4} \cite{Hu8} \cite{HLZ6}-\cite{HLZ8} \cite{Hu14}, 
modular invariance of (logarithmic) intertwining operators \cite{Hu9} \cite{Fio1} \cite{Fio2} \cite{FH}, 
the Verlinde formula \cite{Hu12}, vertex tensor category structures \cite{HL1}-\cite{HL5} \cite{Hu4}
\cite{HLZ1}-\cite{HLZ9} \cite{Hu14} and modular tensor category structures
\cite{Hu13}, have been known to form the foundation of 
conformal field theory and are necessary for a deep understanding of even 
familiar examples such as Wess-Zumino-Witten models, minimal models and orbifold theories obtained from 
such familiar examples. An introductory discussion of some of this work is given in \cite{HL6}.

This 
general theory has been showing its power in the solutions of important problems, including those 
mentioned above and some problems whose solutions have been announced recently.  
We expect that  more problems
will be solved soon using the general theory. On the other hand,
since uniqueness and classification results can only be obtained using the
general theory, the lack of major uniqueness and classification results in this area is an indication that the general 
theory is still far from complete. The author hopes that the problems given here will help to encourage the
further development of the general theory.

For simplicity, by conformal field theories, we shall always mean two-dimensional 
conformal field theories.

\vspace{.7em}
{\bf Acknowledgments}\;\;
The author would like to thank Katrina Barron, Elizabeth Jurisich, Haisheng Li, Antun Milas and 
Kailash Misra 
for organizing this conference and for arranging two problem sessions. The author is grateful 
to Jim Lepowsky for his support starting from the time when the author was still a first year graduate student. 
The author would also like to thank Scott Carnahan, John Duncan, Yasuyuki Kawahigashi, Ching Hung Lam, 
David Radnell, Nils Scheithauer, Eric Schippers,  Wolfgang Staubach 
and Katrin Wendland for conversations and discussions on their works.

\section{The construction of rational conformal field theories satisfying the axioms of Kontsevich-Segal-Moore-Seiberg}

\begin{prob}
Construct rational conformal field theories satisfying the axioms of Kontsevich, G. Segal \cite{Se} and 
Moore-Seiberg \cite{MS1} \cite{MS2}, or at least prove the existence of such rational conformal field theories. In particular, 
construct the Wess-Zumino-Witten models and minimal models or at least prove their existence. 
\end{prob}

In fact, the genus-zero and genus-one parts of the problem have been essentially solved, and in fact, the references
in the introduction include a significant part of the solution of this problem.  For a vertex 
operator algebra satisfying three conditions corresponding to the assumption that the conformal field theory 
is ``rational,'' the operator product expansion, or associativity, of intertwining operators, and the modular invariance
of the spaces of $q$-traces of products of intertwining operators have been proved \cite{Hu8} \cite{Hu9}.
In fact, correlation functions of genus-zero and genus-one chiral rational conformal field theories and their properties 
were given by these results.  The genus-zero and genus-one correlation functions for the full rational conformal field theories
and their properties have also been obtained \cite{HK2} \cite{HK3}. But there is still no 
construction of the higher-genus correlation functions. Also, the state space obtained in these constructions
is only a graded space with a nondegenerate symmetric bilinear or Hermitian form. To construct the 
conformal field theories completely, we also need to construct a locally convex complete topological vector
space with a nondegenerate Hermitian form. 

The construction of the higher-genus correlation functions is the main unsolved problem. If the axioms 
for conformal field theories are assumed, one can see that if a conformal field theory is constructed, then higher-genus 
correlation functions can be expanded as series obtained using intertwining operators. Since we have not constructed
conformal field theories satisfying all the axioms, even though we can still write down these series
using intertwining operators, we cannot use the axioms for conformal field theories to derive the 
convergence of these series. If the convergence of these series can be proved, then the sums of these series
give functions on the Teichm\"{u}ller spaces. Using the associativity of intertwining operators and 
the modular invariance for intertwining operators that have been proved as theorems in  \cite{Hu8} 
and \cite{Hu9}, one can prove that 
these functions on the Teichm\"{u}ller spaces in fact give flat sections of holomorphic vector bundles with flat connections
over the moduli spaces. These are higher-genus chiral correlation functions. Then the construction
of higher-genus full  correlation functions follows trivially. Thus the problem of constructing  the higher-genus correlation functions
is reduced to the problem of proving the convergence of the series above. This convergence problem is the 
generalization in the higher-genus case of the convergence proved in the associativity theorem for intertwining operators in \cite{Hu8} and 
the modular invariance theorem for intertwining operators  
in \cite{Hu9}. 

\begin{prob}
Prove the convergence of these series obtained using intertwining operators to construct higher-genus correlation functions.
\end{prob}

To prove this convergence, we need to prove a conjecture on certain functions on the Teichm\"{u}ller spaces and 
moduli spaces of Riemann
surfaces with parametrized boundaries. Since the conjecture is technical, we shall not state it here. 
In recent years, Radnell, Schippers and Staubach \cite{RS1}-\cite{RS3} \cite{RSW1}-\cite{RSW5}
have made important progress in the study of these Teichm\"{u}ller spaces and 
moduli spaces. They have found the correct class of Riemann surfaces with parametrized boundaries
underlying the definition of conformal field theory by Kontsevich and G. Segal. We expect that the 
conjecture mentioned above on certain functions on the Teichm\"{u}ller spaces and 
moduli spaces will be proved in the near future. Then it can be used to prove the convergence needed in the 
construction of higher-genus correlation functions.

The construction of locally convex complete topological vector
spaces with nondegenerate Hermitian forms is related to the construction of higher-genus correlation functions.
From the axioms of Kontsevich-Segal, we see that if we have indeed constructed a conformal field theory, 
then a Riemann surface $S$ with 
one connected boundary component corresponds   to  a  map from 
the field  of 
 complex numbers to  the state space of
 the  conformal field theory. 
Such a map is equivalent to  an element 
\begin{center}
\includegraphics{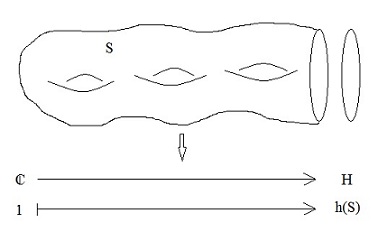}
\end{center}
$h(S)$  in the state
 space. Thus we see that the state space of a conformal field theory
must contain elements corresponding to Riemann surfaces of arbitrary genera. Only after we have a construction of  higher-genus 
correlation functions, we will be able to understand the structure of the state space of the conformal field theory 
completely. The genus-zero full conformal field theory constructed in \cite{HK2} also has a state space. 
But this space is a graded vector space obtained from the tensor products and direct sums of modules for the vertex operator algebra
and does not contain those elements corresponding to Riemann surfaces of arbitrary genera. The author in 
1998 \cite{Hu7.1}  and 2000 \cite{Hu7.2} constructed a locally convex topological completion of a vertex operator
algebra. If we have a construction of  higher-genus 
correlation functions, we can apply the method in \cite{Hu7.1}  and \cite{Hu7.2} to add elements 
corresponding to Riemann surfaces of arbitrary genera to the graded state space of the genus-zero full conformal field theory
in \cite{HK2} to obtain the complete state space of the conformal field theory.

If the graded state space of a genus-zero full conformal field theory has an inner product invariant under the 
conformal fields constructed from intertwining operators (we call it a compatible inner product), then this graded state space also has a Hilbert space 
completion. We have the following conjecture:

\begin{conj}
If the graded state space of a genus-zero full conformal field theory has a compatible  inner product,
then the locally convex topological completion obtained by adding elements corresponding to 
higher-genus Riemann surfaces and the Hilbert space completion are isomorphic as  locally convex topological vector spaces
with grading structures and nondegenerate symmetric Hermitian forms. 
\end{conj}

\section{Cohomology theory for graded vertex algebras and complete reducibility of their modules}

In representation theory, complete reducibility of modules is a basic problem. For an associative algebra,
if the algebra as a module for itself is completely reducible, then all modules for the algebra are completely reducible
and all the irreducible modules appear in the decomposition of the algebra as a direct sum of irreducible modules. 
But for a vertex operator algebra, even if as a module for itself it is irreducible, there might still be many other irreducible 
modules and modules that are not completely reducible. For concrete examples of vertex operator algebras, 
the existing proofs of the complete reducibility of modules have been reduced to the complete reducibility of modules for 
some other algebras or some other properties, not deduced from a general theorem on complete reducibility with easy-to-verify conditions.

An associative algebra is semisimple (equivalent to all modules being completely reducible) if and only if 
its Hochschild cohomology with any bimodule as coefficients is $0$. For a vertex operator algebra,
 the conformal element is irrelevant to the complete reducibility of modules. Thus we need only consider
grading-restricted vertex algebras, which are the same as vertex operator algebras except that 
they do not have conformal elements. On the other hand, the complete reducibility theorem for associative
algebras also applies to commutative associative algebras. There is no need to replace the 
Hochschild cohomology by Harrison cohomology, since the commutativity is irrelevant to the 
complete reducibility of modules. Grading-restricted vertex algebras also have a ``commutativity property,'' which 
is also irrelevant to the complete reducibility of modules. So if there is a similar complete reducibility theorem,
it should be for more general algebras that do not have to satisfy commutativity. In 2012, the author \cite{Hu21}
introduced a notion of meromorphic open-string vertex algebra, which should be viewed as a noncommutative generalization
of the notion of grading-restricted vertex algebra. In 2010, the author \cite{Hu17} \cite{Hu18}  introduced a cohomology theory 
for grading-restricted vertex algebras and proved the basic properties that a cohomology theory must have. 
The construction of this cohomology 
in fact contains two steps. The first step is a construction of a cochain complex that is similar to the Hochschild cochain complex. 
From this cochain complex, we also have a cohomology. The second step is to define a cochain subcomplex of the 
cochain complex above that is similar to the Harrison cochain complex. The cohomology of this cochain subcomplex is 
the cohomology of the grading-restricted vertex algebra. The cohomology constructed in the first step
is in fact the cohomology when we do not consider the commutativity of the grading-restricted vertex algebra,
that is, the cohomology when we view the grading-restricted vertex algebra as a meromorphic open-string vertex algebra.
The cochain complex and cohomology in this step can all be generalized to meromorphic open-string vertex algebras.
Thus for such an algebra, we also have a cohomology theory.

If we regard grading-restricted vertex algebras as analogues of commutative associative algebras, then 
meromorphic open-string vertex algebras are analogues of associative algebras. Using the same analogy and 
what we know about grading-restricted vertex algebras, we have the following conjecture:

\begin{conj}\label{coh-red}
Every grading-restricted generalized module of finite length for a meromorphic open-string vertex algebra
is completely reducible if and only if for $n\in \Z_{+}$, the $n$-th cohomology of the meromorphic open-string vertex algebra
with coefficients in any bimodule is $0$. In particular, every grading-restricted generalized module of finite length for a 
grading-restricted vertex algebra is completely reducible if and only if for $n\in \Z_{+}$, the $n$-th cohomology 
of the grading-restricted vertex algebra viewed as a meromorphic open-string vertex algebra
with coefficients in any bimodule (whose left and right actions can be different)  is $0$.
\end{conj}

In 2015, Qi and the author \cite{HQ} proved, under a convergence assumption,  that if the first cohomology of a 
meromorphic open-string vertex algebra with coefficients in any bimodule is $0$, then
every grading-restricted generalized $V$-module of finite length is completely reducible. Thus the conjecture 
above has been reduced to the convergence assumption and the conjecture 
that if every grading-restricted generalized module of finite length for a meromorphic open-string vertex algebra
is completely reducible, then for $n\in \Z_{+}$, the $n$-th cohomology of the meromorphic open-string vertex algebra
with coefficients in any bimodule is $0$. 

We know that grading-restricted generalized modules of finite length for the vertex operator algebras 
associated to the Wess-Zumino-Witten models, the minimal models, the lattice vertex operator algebras,
and the moonshine module vertex operator algebra are all completely reducible. If Conjecture \ref{coh-red} is true,
then their $n$-th cohomologies for $n\in \Z_{+}$ when they are viewed as meromorphic open-string vertex algebras should be $0$.
Thus to convince ourselves that Conjecture \ref{coh-red} is true, we should first verify that for $n\in \Z_{+}$, 
the $n$-th cohomology of  such a vertex operator algebra as a meromorphic open-string vertex algebra is $0$. 
But the calculation of the cohomology of such a vertex operator algebra is also a nontrivial open problem because
the calculation involves a convergence problem. 

\begin{conj}
Let $V$ be  a vertex operator algebra
associated to a Wess-Zumino-Witten model,  a minimal model, a lattice vertex operator algebra
or the moonshine module vertex operator algebra. Then for $n\in \Z_{+}$, the $n$-th cohomology of $V$
as a meromorphic open-string vertex algebra with coefficients in any $V$-bimodule is $0$. 
\end{conj}

If Conjecture \ref{coh-red} is proved, we would obtain a criterion for the complete reducibility using
the cohomology theory. In addition, the work \cite{HQ} also simplifies this criterion to the 
criterion that the $0$-th cohomology with coefficients in any bimodule is $0$. But this criterion 
must hold for all bimodules. This is still very difficult to verify. Our hope is that this necessary 
and sufficient condition for the complete reducibility can help us to find a criterion that is 
easy to verify.

\begin{prob}
Find a criterion on a grading-restricted vertex algebra itself for the complete reducibility of modules. 
\end{prob}

The author thinks that  the connection between the cohomology and the Killing form of a Lie algebra
might give some hint for the solution of this problem.

\section{The moduli space of conformal field theories}

Moduli spaces always play an important role in mathematics. From the moduli space of Riemann surfaces
to the moduli space of self-dual or anti-self-dual solutions of Yang-Mills equations, many important 
mathematical results have been obtained from the studies of them. The moduli space of conformal field theories 
is also an important mathematical structure. In mathematics, it is closely related to the moduli space
of Calabi-Yau manifolds. In physics, it is closely related to the solution space of 
string theory and it might also be closely related to topological order. But at this moment, it is even not
clear how the topology on this moduli space should be defined. 

\begin{prob}
Study the moduli space of conformal field theories. Give topological and geometric structures on 
the moduli space. Prove that rational conformal field theories are isolated points in the moduli space.
\end{prob}

One of the approaches to studying the moduli space is to develop a deformation theory of 
conformal field theories. A deformation theory will at least help us to understand the topology
of the moduli space. In \cite{Hu18}, the author has proved that the first order deformations of 
a grading-restricted vertex algebra correspond to the second cohomology of this vertex algebra
with coefficients in itself. In fact, the author has also proved that the third cohomology and a convergence 
in this cohomology theory is the obstruction 
for a first order deformation of the vertex algebra to be lifted to a formal deformation. But this proof
involves very complicated and technical calculations. The author hopes that after introducing new concepts
and methods, the calculations can be simplified. So this work has not been published yet. 

We have the following problem:

\begin{prob}
Find the conditions for the formal deformations of grading-restricted vertex algebras to converge
to analytic deformations. Generalize the deformation theory of grading-restricted vertex algebras to 
a deformation theory of genus-zero full conformal field theories. 
\end{prob}

\section{The construction and study of logarithmic conformal field theories}

Though we still need to construct higher-genus correlation functions and prove that the axioms 
of Kontsevich-Segal-Moore-Seiberg are satisfied, we already have a lot of important results on 
rational conformal field theories, including the operator product of intertwining operators \cite{Hu4} \cite{Hu8}, modular 
invariance of intertwining operators \cite{Hu9}, the Verlinde formula \cite{Hu12}, modular tensor category structures
\cite{Hu13}, three-dimensional topological quantum field theories and 
invariants of knots and $3$-manifolds  (combining \cite{Hu13} and \cite{T}), and
genus-zero and genus-one rational chiral and full conformal field theories \cite{HK2} \cite{HK3}. 
But for logarithmic conformal field theories, many of the generalizations of these 
results still have not been proved and in some cases, even the precise formulations of some of the conjectures
are still not known. 

Starting from 2001, Lepowsky, Zhang 
and the author \cite{HLZ1}-\cite{HLZ9} have constructed and developed a tensor category theory 
for logarithmic conformal field theories under suitable natural assumptions and have constructed 
vertex tensor category structures and braided tensor category structures from suitable categories of 
modules for suitable vertex algebras.  
The main results are the proof of the associativity of logarithmic intertwining 
operators (or the logarithmic operator product expansion) and the construction of 
vertex tensor category and braided tensor category structures. Even though they are natural,
the assumptions in this theory are not always trivial to verify. In 2007, the author  \cite{Hu14}
proved that for a vertex operator algebra satisfying the two conditions in the conjecture below, 
the assumptions in the work \cite{HLZ1}-\cite{HLZ9}
 are satisfied and thus in this case, the category of grading-restricted generalized modules for 
the vertex operator algebra has natural structures of a vertex tensor category and of a braided tensor category. 

In 2009, the author \cite{Hu15} proposed the following conjecture:

\begin{conj}\label{log-rigidity}
Assume that $V$ is a simple vertex operator algebra 
satisfying the following conditions:
\begin{enumerate}
\item $V_{(0)}=\C\one$, $V_{(n)}=0$ 
for $n<0$ and the contragredient $V'$, as a $V$-module, is equivalent to $V$.

\item $V$ is $C_{2}$-cofinite, that is, $\dim V/C_{2}(V)<\infty$, where $C_{2}(V)$ is the subspace of 
$V$ spanned by the elements of the form $\res_{x}x^{-2}Y(u, x)v$ for $u, v\in V$ and $Y: V\otimes V\to V[[x, x^{-1}]]$
is the vertex operator map for $V$.
\end{enumerate}
Then the braided tensor category given in \cite{Hu14} based on the work \cite{HLZ1}-\cite{HLZ9}
is rigid. 
\end{conj}

In the case of rational conformal field theories, the proof of the rigidity \cite{Hu13} used the modular invariance \cite{Hu9}. 
For a $C_{2}$-cofinite vertex operator algebra without elements of negative weight, the work of 
Fiordalisi \cite{Fio1} \cite{Fio2} and a joint paper in preparation by Fiordalisi and the author \cite{FH}
have proved that modular invariance for logarithmic intertwining operators also holds. We expect that 
this modular invariance can be used to prove the rigidity conjecture above. 

For a rational conformal field theory, the corresponding category of modules for the vertex operator algebra
is not only a rigid braided tensor category, but also satisfies a nondegeneracy property, so that together with 
some other minor properties that it has, it is a 
modular tensor category. But this nondegeneracy property is formulated using irreducible objects 
(irreducible modules)  in the category, since for a rational conformal field theory, the category is
semisimple, that is, every object is a direct sum of irreducible objects. For a logarithmic conformal 
field theory, since the category is not semisimple, it is not sufficient to consider only irreducible modules. 
So giving the correct definition of modular tensor category in the case of logarithmic conformal field theories
is also an important unsolved problem. 

\begin{prob}
Give a definition of modular tensor category in the nonsemisimple case and prove that the braided 
tensor category in Conjecture \ref{log-rigidity} has such a modular tensor category structure. From
such a nonsemisimple modular tensor category, can we still construct a $3$-dimensional topological 
quantum field theory and invariants of knots and $3$-manifolds?
\end{prob}

The associativity of logarithmic intertwining operators and the modular invariance of logarithmic intertwining 
operators discussed above gave the genus-zero and genus-one chiral correlation functions in the corresponding 
logarithmic conformal field theory. But we do not know how to construct the genus-zero and genus-one full correlation
functions for the logarithmic full conformal field theory. 

\begin{prob}
For a vertex operator algebra satisfying the conditions in Conjecture \ref{log-rigidity}, 
is it possible to construct genus-zero and genus-one full correlation
functions from the genus-zero and genus-one chiral correlation functions? From these 
genus-zero and genus-one correlation functions, is it possible to construct a conformal field theory 
satisfying the axioms of Kontsevich-Segal, except for the axioms involving unitarity?
\end{prob}

\section{Orbifold conformal field theories}

Orbifold conformal field theories play an important role in the construction and application of conformal field theory. 
The moonshine module vertex operator algebra constructed by Frenkel, Lepowsky and Meurman \cite{FLM}
is in fact the first example of orbifold conformal field theories. To study orbifold conformal field theories, we need 
to first study twisted modules for vertex operator algebras. But the construction and study of twisted modules 
are still far away from the complete construction of orbifold conformal field theories. First, the author has the 
following conjecture:

\begin{conj}\label{twisted-int-op}
Let $V$ be a vertex operator algebra satisfying the two conditions in Conjecture \ref{log-rigidity}. Assume that 
in addition every grading-restricted generalized $V$-module is completely reducible. Let $G$ be a finite group 
of automorphisms of $V$. Then the twisted intertwining operators among the $g$-twisted $V$-modules for all $g\in G$
satisfy the associativity, commutativity and modular invariance properties.
\end{conj}

If we replace the condition requiring that the vertex operator algebra is rational in a conjecture on the 
category of twisted modules (see Example 5.5 in \cite{Ki}) by the three conditions in Conjecture \ref{twisted-int-op},
then we obtain the following conjecture:

\begin{conj}
Let $V$ be a vertex operator satisfying the three conditions in Conjecture \ref{twisted-int-op}
and let $G$ be a finite group of automorphisms of $V$. The the category of $g$-twisted $V$-modules for 
all $g\in G$ is a $G$-crossed (tensor) category in the sense of Turaev \cite{T2}. 
\end{conj}

These two conjectures are both under the complete reducibility assumption and are also about 
finite groups of automorphisms of $V$. In the case that grading-restricted generalized $V$-modules 
are not complete reducible and $G$ is not finite, we have the following conjecture and problem:

\begin{conj}
Let $V$ be a vertex operator algebra satisfying the two conditions in Conjecture \ref{log-rigidity}
and let $G$ be a finite group of automorphisms of $V$. Then the twisted logarithmic intertwining operators 
among the $g$-twisted $V$-modules for all $g\in G$
satisfy the associativity, commutativity and modular invariance properties.
\end{conj}

\begin{prob}
Let $V$ be  a vertex operator algebra and let $G$ be a group of automorphisms of $V$. If $G$ is an infinite group,
under what conditions do the twisted logarithmic intertwining operators 
among the $g$-twisted $V$-modules for all $g\in G$
satisfy the associativity, commutativity and modular invariance properties? Under what conditions is the category of 
$g$-twisted $V$-modules for all
$g\in G$  a $G$-crossed (tensor) category?
\end{prob}

\section{The uniqueness of the moonshine module vertex operator algebra and the classification of 
meromorphic rational conformal field theories of central charge $24$}

An important problem in the study of conformal field theories is the classification of rational conformal 
field theories. But for general rational conformal field theories, since this problem might be related to 
the classification of finite simple groups, it is expected to be very difficult. A relatively practical problem is the classification of 
meromorphic rational conformal field theories. Here by a meromorphic rational conformal field theory we mean 
a rational conformal field theory whose chiral vertex operator algebra is the only irreducible module for itself. 
The classification of meromorphic rational conformal field theories is equivalent to the classification of 
the corresponding chiral vertex operator algebras. 

Meromorphic rational conformal field theories of central charge 24 are especially important. In 1988, 
Frenkel, Lepowsky and Meurman proposed the following well-known uniqueness conjecture for the 
moonshine module vertex operator algebra:

\begin{conj}\label{uniqueness-moonshine}
Let $V$ be a vertex operator algebra satisfying the following three conditions:
\begin{enumerate}
\item The central charge of $V$ is $24$.

\item $V$ has no nonzero elements of weight $1$. 

\item Any irreducible $V$ module is equivalent to $V$ as a $V$-module and every $V$-module is 
completely reducible.
\end{enumerate}
Then $V$ must be isomorphic to the moonshine module vertex operator algebra. 
\end{conj}

If Conjecture \ref{uniqueness-moonshine} is true, then the Monster group can be defined as the automorphism 
group of a vertex operator algebra satisfying the three conditions in this conjecture.

In 1992, based on the study and conjectures on the Lie algebras obtained from 
the homogeneous subspaces of weight $1$ of meromorphic rational conformal field theories, 
Schellekens \cite{Sc} proposed a classification conjecture of meromorphic rational conformal field theories 
with nonzero homogeneous subspaces of weight $1$ and of central charge $24$. The conjecture says that there are 
$70$ such  meromorphic rational conformal field theories. Together with the moonshine module vertex operator algebra,
conjecturally there is a total of $71$ meromorphic rational conformal field theories of central charge $24$.
When Schellekens proposed the conjectures, there were only $39$ vertex operator algebras for such meromorphic rational conformal field theories
constructed. Also, until 2016, even the statement that there are only $70$ possible nonzero Lie algebras 
obtained from the homogeneous subspaces of weight $1$ of meromorphic rational conformal field theories of central charge $24$
was a conjecture. 

In recent years,  important progress has been made on the existence part of Schellekens' conjecture. 
Lam \cite{La}, Lam-Shimakura \cite{LS1} \cite{LS2} \cite{LS3},  Miyamoto \cite{Miy3}, Shimakura-Sagaki \cite{SS}, 
van Ekeren-M\"{o}ller-Scheithauer \cite{EMS} have obtained $68$ of the $70$ 
vertex operator algebras of central charge $24$ with nonzero homogeneous subspaces of weight $1$
in Schellekens' classification conjecture. There are still two of them that need to be constructed. 
In \cite{EMS}, van Ekeren, M\"{o}ller and Scheithauer also proved that 
the Lie algebra obtained from the homogeneous subspace of weight $1$ of 
such a vertex operator algebra must be one of the $70$ Lie algebras in the list of Schellekens. 

On the other hand, there is no real progress towards the uniqueness part of the Schellekens' classification conjecture, except for the 
uniqueness, proved in \cite{EMS}, of the Lie algebras obtained from these vertex operator algebras.
Here we state this uniqueness part as the main remaining conjecture in this classification problem. 

\begin{conj}\label{uniqueness}
Let $V$ be a vertex operator algebra satisfying the first and the third conditions in Conjecture \ref{uniqueness-moonshine}.
Then $V$ is determined uniquely by the Lie algebra structure on its homogeneous subspace of weight $1$.
\end{conj}

These two conjectures might need stronger complete reducibility conditions. It is very likely that we might 
need to assume all weak modules are completely reducible.

In the case that the homogeneous subspace of weight $1$ of such a vertex operator algebra is not $0$, we can 
use all the possible vertex operator
algebras generated by the corresponding Lie algebras in the list of Schellekens to study Conjecture \ref{uniqueness}.
The difficulty of the uniqueness of the moonshine module vertex operator algebra is that there is no such Lie algebra
that one can use. The author believes that the proofs of these conjectures need more powerful theories and sophisticated tools 
of vertex operator algebras that have been developed and are being developed.

\section{Calabi-Yau superconformal field theories}

In 1985,  physicists, including Friedan,  Candelas, Horowitz, Strominger, Witten,
Alvarez-Gaum\'{e},  Coleman, Ginsparg
 \cite{F} \cite{CHSW} \cite{ACG}  (see also the paper \cite{NS} by Nemeschansky and  Sen),  proposed a conjecture that 
the quantum nonlinear $\sigma$-model with a Calabi-Yau manifold as the target space is an 
$N=2$ superconformal field theory, that is, a conformal field theory with $N=2$ superconformal structures. 

In 1987, Gepner \cite{Ge} obtained an $N=2$ superconformal field theory from $N=2$ superconformal 
minimal models (now called the Gepner model). He also conjectured that this $N=2$ superconformal field theory
should be isomorphic to the quantum nonlinear $\sigma$-model with the Fermat quintic threefold 
in four-dimensional complex projective space as the target space. 

To study these conjectures using the representation theory of vertex operator algebras, the first step 
is to construct the corresponding $N=2$ vertex operator superalgebras, that is, vertex operator algebras with 
$N=2$ superconformal structures. 

\begin{prob}\label{CY-funct}
Construct  a functor from the category of Calabi-Yau manifolds to the category of $N=2$ vertex operator superalgebras
such that the Fermat quintic threefold 
corresponds to the $N=2$ vertex operator superalgebra for the Gepner model and
such that deformations of Calabi-Yau manifolds correspond
to deformations of $N=2$ superconformal field theories. 
\end{prob}

Although there have been constructions of $N=2$ vertex operator superalgebras from some Calabi-Yau manifolds,
there is no general construction that gives the functor in the problem above.

In \cite{H13}, the author constructed a sheaf $\mathcal{V}$ of meromorphic open-string vertex algebras 
on a Riemannian manifold and a sheaf of left modules for $\mathcal{V}$ generated by the sheaf of smooth functions.
The Laplacian on the Riemannian manifold was shown in \cite{H13} to be a component of a vertex operator
acting on the space of smooth functions. In particular, we have sheaves of left modules for $\mathcal{V}$
generated by eigenfunctions 
for the Laplacian. Since the state space of the nonlinear $\sigma$-models with a Riemannian manifold as the target
must contain eigenfunctions for the Laplacian, this construction must be related to the nonlinear $\sigma$-model.
Certainly for a Riemannian manifold,
the corresponding nonlinear $\sigma$-model is in general not a conformal field theory.
But this corresponds exactly to the fact that in the general case, we obtain sheaves of left modules for $\mathcal{V}$
generated by eigenfunctions instead of modules for a sheaf of 
vertex operator algebras. We  expect that when the Riemannian manifold is a 
Calabi-Yau manifold, we shall be able to construct sheaves of modules generated by eigenfunctions,
eigenforms and eigenspinors for a sheaf of $N=2$ superconformal vertex operator superalgebra,
not just for $\mathcal{V}$. We hope that this construction will
eventually lead to a construction of the functor discussed in Problem \ref{CY-funct}.

\begin{prob}
If $N=2$ vertex operator superalgebras corresponding to Calabi-Yau manifolds can be constructed, study the
representation theory of these $N=2$ vertex operator superalgebras. Use this representation theory to 
construct the corresponding $N=2$ superconformal field theories, including the proof of the operator product expansion of 
intertwining operators, the proof of modular invariance, the construction of modular tensor categories (in a suitable sense),
and the construction of full $N=2$ superconformal field theories.
\end{prob}

For Calabi-Yau manifolds, the construction of full $N=2$ superconformal field theories is very important because
many conjectures given by physicists (for example, quantum cohomology and mirror symmetry) are obtained from
full $N=2$ superconformal field theories but cannot be obtained from only chiral $N=2$ superconformal field theories. 

\begin{prob}
Suppose that $N=2$ vertex operator superalgebras corresponding to Calabi-Yau manidfolds can be constructed, that
the basic results in the representation theory of these algebras can be obtained and that the full $N=2$ superconformal field theories
can be constructed. Formulate and prove the conjectures (for example, quantum cohomology and mirror symmetry) by physicists
on Calabi-Yau manifolds by turning physicists' intuition into mathematical methods using these results and constructions.
\end{prob}

Recently, there have been some very interesting developments in the study of the $N=2$ vertex operator superalgebras
corresponding to $K3$ surfaces. Early in 2001,
Wendland \cite{We} studied the full $N=2$ superconformal field theory corresponding to a special $K3$ surface. 
She found that it has a very large but finite automorphism group. In 2013, 
Gaberdiel, Taormina, Volpato and Wendland \cite{GTVW} used the different descriptions of this particular 
full $N=2$ superconformal field theory to determine completely its huge finite automorphism group. 
In 2015, Duncan and Mack-Crane \cite{DM2} discovered that the state space of the Neveu-Schwarz sector of this 
full $N=2$ superconformal field theory as a module for the Virasoro algebra is equivalent to 
the moonshine module vertex operator superalgebra $V^{S\natural}$ for the Conway group that they had 
constructed in 2014 \cite{DM1}. They also proved that the state space of the Ramond sector of this 
full $N=2$ superconformal field theory as a module for the Virasoro algebra is equivalent to a twisted $V^{S\natural}$-module.
Using the connections established by these works, Taormina and Wendland have gone further to use 
the vertex operator superalgebra $V^{S\natural}$ to describe the full $N=2$ superconformal field theory corresponding to 
this special $K3$ surface. We hope that these studies will provide important mathematical ideas and examples for 
the future construction of full $N=2$ superconformal field theories corresponding to Calabi-Yau manifolds,
including hints and ideas for the solutions of the problems above.

\section{The relation between the approaches of vertex operator algebras and conformal nets}

There are many different methods for the study of conformal field theories. But they can all be 
classified as belonging to one of two types of methods. The first type is the method of 
the representation theory of vertex operator algebras. The other is the method of 
conformal nets. The problems discussed above are all problems using the first type of method. 
For the method of conformal nets, the author recommends the survey \cite{Ka} by Kawahigashi. 
For some results on the connection between the two methods, see \cite{CKLW} by 
Carpi, Kawahigashi,  Longo and Weiner. For a functional-analytic study of vertex operator 
algebras and their modules, see \cite{Hu7.1} and \cite{Hu7.2}. Both methods gave
modular tensor category structures corresponding to rational conformal field theories. 
For the study of other problems,  they have different advantages.

\begin{prob}
Find  a direct connection between the method of the representation theory of vertex operator algebras
and the method of conformal nets. Prove that at least for rational conformal field theories, they are 
equivalent. 
\end{prob}


\bibliographystyle{amsalpha}

\begin{thebibliography}{KWAK2}

\bibitem[ACG]{ACG}
L. Alvarez-Gaum\'{e}, S. Coleman and P. Ginsparg,
Finiteness of Ricci flat $N=2$ supersymmetric $\sigma$-models, {\it Comm. Math. Phys.}
{\bf 103} (1986), 423--430.

\bibitem[CHSW]{CHSW}
P. Candelas, G. Horowitz, A. Strominger and E. Witten,
Vacuum configurations for superstrings, {\it Nucl. Phys.} {\bf B258} (1985), 46--74.

\bibitem[CKLW]{CKLW}
S. Carpi, Y. Kawahigashi, R. Longo, M. Weiner,
From vertex operator algebras to conformal nets and back,
{\it Memoirs Amer. Math. Soc.}, to appear; arXiv:1503.01260.

\bibitem[DM-C1]{DM1}
J. F. R. Duncan and S. Mack-Crane, Derived equivalences of $K3$
surfaces and twined elliptic genera, {\it Res. Math. Sci.} {\bf 3}
(2016), 3:1.

\bibitem[DM-C2]{DM2}
J. F. R. Duncan and S. Mack-Crane, The moonshine module for Conway's group,
{\it Forum Math. Sigma} {\bf 3} (2015), e10.

\bibitem[EMS]{EMS}
J. van Ekeren, S. M\"{o}ller, N. R. Scheithauer,
Construction and classification of holomorphic vertex operator algebras,
to appear; arXiv:1507.08142

\bibitem[Fi1]{Fio1}
F. Fiordalisi, Logarithmic intertwining operators and genus-one
correlation functions, Ph.D. thesis, Rutgers University, 2015.

\bibitem[Fi2]{Fio2}
F. Fiordalisi, Logarithmic intertwining operators and genus-one
correlation functions, {\it Comm. Contemp. Math}, to appear; arXiv:1602.03250.

\bibitem[FH]{FH}
F. Fiordalisi and Y.-Z. Huang, Modular invariance for
logarithmic intertwining operators, in preparation.

\bibitem[FLM]{FLM}
I.~B. Frenkel, J.~Lepowsky and A.~Meurman,
{\em Vertex Operator Algebras and the Monster},
Pure and Appl. Math., Vol. 134,  Academic Press,  Boston, 1988.

\bibitem[Fr]{F}
D. Friedan, Nonlinear Models in $2+\varepsilon$ Dimensions,
{\it Annals of Physics} {\bf 163} (1985), 318--419.

\bibitem[GTVW]{GTVW}
M. Gaberdiel, A. Taormina, R. Volpato, K.
Wendland, A $K3$ sigma model with $\Z_2^8$: $\mathbb{M}_20$ symmetry,
{\it JHEP}  {\bf 2}
(2014), 022.

\bibitem[G]{Ge}
D. Gepner, Space-time supersymmetry in compactified string theory and
superconformal models, {\it Nucl. Phys.} {\bf B296} (1988), 757--778.

\bibitem[H1]{Hu4}
Y.-Z. Huang, A theory of tensor products for module categories for a
vertex operator algebra, IV, {\em J. Pure Appl. Alg.} 100 (1995)
173--216.

\bibitem[H2]{Hu7.1}
Y.-Z. Huang, A functional-analytic theory of vertex (operator) algebras, I,
{\it Comm. Math. Phys.} {\bf 204} (1999), 61--84.

\bibitem[H3]{Hu7.2}
Y.-Z. Huang, A functional-analytic theory of vertex (operator) algebras, II,
{\it Comm. Math. Phys.} {\bf 242} (2003), 425--444.

\bibitem[H4]{Hu8} Y.-Z. Huang, Differential equations and
intertwining operators, {\em Comm. Contemp. Math.} {\bf 7} (2005),
375--400.

\bibitem[H5]{Hu9}
Y.-Z. Huang,  Differential equations, duality and modular invariance,
{\em Comm. Contemp. Math.} {\bf 7} (2005), 649--706.

\bibitem[H6]{Hu12} Y.-Z. Huang, Vertex operator
algebras and the Verlinde conjecture, {\em Comm. Contemp. Math.}
{\bf 10} (2008), 103--154.

\bibitem[H7]{Hu13} Y.-Z. Huang, Rigidity and modularity of vertex
tensor categories, {\em Comm. Contemp. Math.} {\bf 10} (2008), 871--911.

\bibitem[H8]{Hu14} Y.-Z. Huang, Cofiniteness conditions, projective
covers and the logarithmic tensor product theory, {\em J. Pure
Appl. Alg.} {\bf 213} (2009), 458--475.

\bibitem[H9]{Hu15}
Y.-Z. Huang, Representations of vertex operator algebras and braided
finite tensor categories, in: {\it Vertex Operator Algebras and Related
Topics, An International Conference in Honor of Geoffery Mason's 60th Birthday},
ed. M. Bergvelt, G. Yamskulna and W. Zhao, Contemporary Math., Vol. 497,
Amer. Math. Soc., Providence, 2009, 97--111.

\bibitem[H10]{Hu17}
Y.-Z. Huang, A cohomology theory of grading-restricted vertex algebras,
{\it Comm. Math. Phys.} {\bf 327} (2014), 279--307.


\bibitem[H11]{Hu18}
Y.-Z. Huang, First and second cohomologies of grading-restricted vertex algebras,
{\it Comm. Math. Phys.} {\bf 327} (2014), 261--278.

\bibitem[H12]{Hu21}
Y.-Z. Huang, Meromorphic open-string vertex algebras, {\it J. Math. Phys.} {\bf 54}
(2013), 051702.

\bibitem[H13]{H13}
Y.-Z. Huang, Meromorphic open-string vertex algebras and Riemannian manifolds, 
to appear; arXiv:1205.2977. 

\bibitem[HK1]{HK2}
Y.-Z. Huang and L. Kong, Full field algebras, {\it Comm. Math. Phys.} {\bf 272}
(2007), 345--396.

\bibitem[HK2]{HK3}
Y.-Z. Huang and L. Kong, Modular invariance for conformal full field algebras,
{\it Trans. Amer. Math. Soc.} {\bf 362} (2010), 3027--3067.

\bibitem[HL1]{HL1}
Y.-Z. Huang and J. Lepowsky, Toward a theory of tensor product for representations
of a vertex operator algebra, in {\it Proc. 20th Intl. Conference on
Diff. Geom. Methods in Theoretical Physics, New York, 1991}, ed. S. Catto and
A. Rocha, World Scientific, Singapore, 1992, Vol. 1, 344--354.

\bibitem[HL2]{HL2}
Y.-Z. Huang and J. Lepowsky, Tensor products of modules for a vertex
operator algebras and vertex tensor categories, in:
     {\em Lie Theory and Geometry,
in honor of Bertram Kostant,}
ed. R. Brylinski, J.-L. Brylinski, V. Guillemin, V. Kac,
Birkh\"{a}user, Boston, 1994, 349--383.

\bibitem[HL3]{HL3}
Y.-Z. Huang and J. Lepowsky, A theory of tensor products for module
categories for a vertex operator algebra, I, {\em Selecta Mathematica
(New Series)} {\bf 1} (1995), 699--756.

\bibitem[HL4]{HL4}
Y.-Z. Huang and J. Lepowsky, A theory of tensor products for module
categories for a vertex operator algebra, II, {\em Selecta Mathematica
(New Series)} {\bf 1} (1995), 757--786.

\bibitem[HL5]{HL5}
Y.-Z. Huang and J. Lepowsky, A theory of tensor
products for module categories for a vertex operator algebra, III,
{\em J. Pure Appl. Alg.} {\bf 100} (1995) 141--171.

\bibitem[HL6]{HL6}
Y.-Z. Huang and J. Lepowsky, Tensor categories and the mathematics of rational 
and logarithmic conformal field theory, {\it  J. Phys.}  {\bf A46} (2013), 494009. 

\bibitem[HLZ1]{HLZ1} Y.-Z.~Huang, J.~Lepowsky and L.~Zhang, A
logarithmic generalization of tensor product theory for modules for a
vertex operator algebra, {\em Internat. J. Math.} {\bf 17} (2006),
975--1012.

\bibitem[HLZ2]{HLZ2}
Y.-Z.~Huang, J.~Lepowsky and L.~Zhang,
Logarithmic tensor category theory for generalized modules for a conformal vertex 
algebra, I: Introduction and strongly graded algebras and their generalized modules,
in: Conformal Field Theories and Tensor Categories, Proceedings of a Workshop Held 
at Beijing International Center for Mathematics Research, ed. C. Bai, J. Fuchs, Y.-Z. 
Huang, L. Kong, I. Runkel and C. Schweigert, Mathematical Lectures from Beijing 
University, Vol. 2, Springer, New York, 2014, 169--248.

\bibitem[HLZ3]{HLZ3} Y.-Z.~Huang, J.~Lepowsky and L.~Zhang, Logarithmic
tensor category theory, II: Logarithmic formal calculus
and properties of logarithmic intertwining operators, to appear; arXiv:1012.4196.

\bibitem[HLZ4]{HLZ4} Y.-Z.~Huang, J.~Lepowsky and L.~Zhang, Logarithmic
tensor category theory, III: Intertwining maps and tensor
product bifunctors, to appear; arXiv:1012.4197.

\bibitem[HLZ5]{HLZ5} Y.-Z.~Huang, J.~Lepowsky and L.~Zhang, Logarithmic
tensor category theory, IV: Constructions of tensor
product bifunctors and the compatibility conditions, to appear; arXiv:1012.4198.

\bibitem[HLZ6]{HLZ6} Y.-Z.~Huang, J.~Lepowsky and L.~Zhang, Logarithmic
tensor category theory, V: Convergence condition for
intertwining maps and the corresponding compatibility
condition, to appear; arXiv:1012.4199.

\bibitem[HLZ7]{HLZ7} Y.-Z.~Huang, J.~Lepowsky and L.~Zhang, Logarithmic
tensor category theory, VI: Expansion condition, associativity of logarithmic
intertwining operators, and the associativity isomorphisms, to appear; arXiv:1012.4202.

\bibitem[HLZ8]{HLZ8} Y.-Z.~Huang, J.~Lepowsky and L.~Zhang, Logarithmic
tensor category theory, VII: Convergence and extension
properties and applications to expansion for intertwining
maps, to appear; arXiv:1110.1929.

\bibitem[HLZ9]{HLZ9} Y.-Z.~Huang, J.~Lepowsky and L.~Zhang, Logarithmic
tensor category theory, VIII: Braided tensor category
structure on categories of generalized modules for a
conformal vertex algebra, to appear; arXiv:1110.1931.

\bibitem[HQ]{HQ}
Y.-Z. Huang and F. Qi, The first cohomology, derivations and the reductivity of a (meromorphic open-string)
vertex algebra, to appear.

\bibitem[Ka]{Ka}
Y. Kawahigashi, Conformal Field Theory, Tensor Categories and Operator Algebras
{\it J. Phys.} {\bf A48} (2015), 303001.

\bibitem[Ki]{Ki}
A. Kirillov, Jr,
On G--equivariant modular categories, arXiv:math/0401119.

\bibitem[L]{La}
C.H. Lam, On the constructions of holomorphic vertex operator algebras of
central charge $24$, {\it Comm. Math. Phys.} {\bf 305} (2011), 153–-198.

\bibitem[LS1]{LS1}
C.H. Lam and H. Shimakura， Quadratic spaces and holomorphic framed vertex
operator algebras of central charge 24, {\it Proc. Lond. Math. Soc.} {\bf 104}
(2012), 540-–576.

\bibitem[LS2]{LS2}
C.H. Lam and H. Shimakura, Classification of holomorphic framed vertex operator
algebras of central charge $24$, {\it Amer. J. Math.} {\bf 137} (2015), 111–-137.

\bibitem[LS3]{LS3}
C.H. Lam and H. Shimakura，Orbifold construction of holomorphic vertex
operator algebras associated to inner automorphisms,  {\it Comm. Math. Phys.},
to appear; arXiv:1501.05094.

\bibitem[M]{Miy3}
M. Miyamoto, A $\mathbb{Z}_{3}$-orbifold theory of lattice vertex operator algebra
and $\mathbb{Z}_{3}$-orbifold constructions,
in {\it Symmetries, integrable systems and representations}, Springer Proc.
Math. Stat., Vol. 40, Springer, Heidelberg, 2013, 319–-344.

\bibitem[MS1]{MS1}
G.~Moore and N.~Seiberg,
Polynomial equations for rational conformal field theories,
{\em Phys. Lett.} {\bf B212} (1988), 451--460.

\bibitem[MS2]{MS2}
G.~Moore and N.~Seiberg,
Classical and quantum conformal field theory,
{\em Comm. Math. Phys.} {\bf 123} (1989), 177--254.

\bibitem[NS]{NS}
D. Nemeschansky and A. Sen,
Conformal invariance of supersymmetric $\sigma$-models on Calabi-Yau manifolds,
{\it Phys. Lett.} {\bf B178} (1986), 365--369.

\bibitem[RS1]{RS1}
D. Radnell and E. Schippers, Quasisymmetric sewing in rigged Teichm\"{u}ller space,
{\it Comm. Contemp. Math.} {\bf 8} (2006), 481--534.

\bibitem[RS2]{RS2}
D. Radnell and E. Schippers, A complex structure on the set of quasiconformally extendible nonoverlapping
mappings into a Riemann surface, {\it J. Anal. Math.} {\bf 108} (2009), 277--291.

\bibitem[RS3]{RS3}
D. Radnell and E. Schippers, Fiber structure and local coordinates for the
Teichm\"{u}ller space of a bordered Riemann surface, {\it Conform. Geom. Dyn.} {\bf 14} (2010), 14--34.

\bibitem[RSS1]{RSW1}
D. Radnell, E. Schippers and W. Staubach,
A Hilbert manifold structure on the Weil-Petersson class
Teichm\"{u}ller space of bordered Riemann surfaces,
{\it Comm. Contemp. Math.} {\bf 17}(2015), 1550016.

\bibitem[RSS2]{RSW2}
D. Radnell, E. Schippers and W. Staubach,
Weil-Petersson class non-overlapping mappings into a
Riemann surface, {\it Comm. Contemp. Math.}, {\bf 18} (2016), 1550060.


\bibitem[RSS3]{RSW3}
D. Radnell, E. Schippers and W. Staubach,
Quasiconformal maps of bordered Riemann surfaces with
$L^2$ Beltrami differentials, {\it J. Anal. Math.},  to appear.

\bibitem[RSS4]{RSW4}
D. Radnell,  E. Schippers and W. Staubach,
Convergence of the Weil-Petersson metric on the Teichm\"{u}ller space
of Bordered Riemann Surfaces, {\it Comm. Contemp. Math.}, to appear.

\bibitem[RSS5]{RSW5}
D. Radnell,  E. Schippers and W. Staubach,
Quasiconformal Teichm\"{u}ller theory as an analytical
foundation for two dimensional conformal field theory, in; {\it Proceedings of the
Conference on Lie Algebras, Vertex Operator Algebras,
and Related Topics, held at University of Notre Dame, Notre Dame, Indiana,
August 14-18, 2015}, ed. K. Barron, E. Jurisich, H. Li, A. Milas, K. C. Misra,
Contemp. Math, American Mathematical Society, Providence,
RI, to appear.

\bibitem[SS]{SS}
D. Sagaki and H. Shimakura, Application of a Z3-orbifold construction to the lattice vertex operator
algebras associated to Niemeier lattices, {\it Trans. Amer. Math. Soc.}
{\bf 368} (2016), 1621-1646.

\bibitem[Sc]{Sc}
A. N. Schellekens， Meromorphic $c=24$ conformal field theories,
{\it Comm. Math. Phys.} {\bf 153} (1993) 159--186.

\bibitem[Se]{Se}
G. Segal, The definition of conformal field theory, in: {\it Topology, Geometry and 
Quantum Field Theory: Proceedings of the 2002 Oxford Symposium in Honour of 
the 60th Birthday of Graeme Segal}, ed. U. Tillmann, London Mathematical Society 
Lecture Note Series, Vol. 308, Cambridge University Press, Cambridge, 2004, 421--577.

\bibitem[T1]{T}
V.  Turaev, {\em Quantum Invariants of Knots and $3$-manifolds},
de Gruyter Studies in Math., Vol. 18,
Walter de Gruyter, Berlin, 1994.

\bibitem[T2]{T2}
V. Turaev, Homotopy field theory in dimension 3 and crossed group-categories,
arxiv:math.GT/0005291.

\bibitem[W]{We}
K. Wendland, Orbifold constructions of $K3$: A link between
conformal field theory and geometry, in: {\it Proceedings of the
Conference on Mathematical Aspects of Orbifold String Theory held at
the University of Wisconsin, Madison, WI, May 4–8, 2001},  ed. A. Adem, J. Morava
and Y. Ruan, Contemp. Math., Vol. 310, 2002, 333-358.

\end{thebibliography}

\end{document}